\documentclass{article}
\usepackage{color}

\catcode`\@=11 \@addtoreset{equation}{section}

\catcode`\@=12

\newtheorem{Theorem}{Theorem}[section]
\newtheorem{Proposition}{Proposition}[section]
\newtheorem{Lemma}{Lemma}[section]
\newtheorem{Corollary}{Corollary}[section]

\newcommand{\bTheorem}[1]{
\begin{Theorem} \label{T#1} }
\newcommand{\eT}{\end{Theorem}}

\newcommand{\bProposition}[1]{
\begin{Proposition} \label{P#1}}
\newcommand{\eP}{\end{Proposition}}

\newcommand{\bLemma}[1]{
\begin{Lemma} \label{L#1} }
\newcommand{\eL}{\end{Lemma}}

\newcommand{\bCorollary}[1]{
\begin{Corollary} \label{C#1} }
\newcommand{\eC}{\end{Corollary}}

\newcommand{\bFormula}[1]{
\begin{equation} \label{#1}}
\newcommand{\eF}{\end{equation}}

\newcommand{\Ov}[1]{\overline{#1}}

\newcommand{\vr}{\varrho}

\newcommand{\vt}{\vartheta}
\newcommand{\vu}{\vc{u}}
\newcommand{\vc}[1]{{\bf #1}}

\newcommand{\Div}{{\rm div}_x}
\newcommand{\Grad}{\nabla_x}

\newcommand{\tn}[1]{\mbox {\F #1}}
\newcommand{\dx}{{\rm d} {x}}
\newcommand{\dt}{{\rm d} t }

\newcommand{\intO}[1]{\int_{\Omega} #1 \ \dx}

\font\F=msbm10 scaled 1000
\newcommand{\R}{\mbox{\F R}}

\date{}
\makeindex
\begin{document}


\title{Weak-strong uniqueness property for the full Navier-Stokes-Fourier system}
\author{Eduard Feireisl \thanks{The work was supported by Grant 201/09/
0917 of GA \v CR as a part of the general research programme of the
Academy of Sciences of the Czech Republic, Institutional Research
Plan AV0Z10190503.}\and Anton{\' \i}n
Novotn\' y }

\maketitle

\bigskip

\centerline{Institute of Mathematics of the Academy of Sciences of
the Czech Republic} \centerline{\v Zitn\' a 25, 115 67 Praha 1,
Czech Republic}

\medskip

\centerline{IMATH, Universit\' e du Sud Toulon-Var, BP 132, 839 57
La Garde, France}

\medskip

\begin{abstract}
The Navier-Stokes-Fourier system describing the motion of a compressible, viscous and heat conducting fluid is known to possess global-in-time weak solutions for any initial data of finite energy. We show that a weak solution coincides with the strong solution, emanating from the same initial data, as long as the latter exists. In particular, strong solutions are unique within the class of weak solutions.
\end{abstract}

\bigskip \noindent
\section{Introduction}
\label{i}

The concept of \emph{weak solution} in fluid dynamics was
introduced in the seminal paper by Leray \cite{LER} in the context
of incompressible, linearly viscous fluids. Despite a concerted
effort of generations of mathematicians, the weak solutions still
represent the only available framework for studying the time
evolution of fluid systems subject to large data on long time
intervals, see Fefferman \cite{Feff}. The original ideas of Leray
have been put into the elegant framework of generalized
derivatives (distributions) and the associated abstract function
spaces of Sobolev type, see Ladyzhenskaya \cite{LAD}, Temam
\cite{TEM}, among many others. More recently, Lions \cite{LI4}
extended the theory to the class of barotropic flows (see also
\cite{EF70}), and, finally, the existence of global-in-time
generalized solutions for the full Navier-Stokes-Fourier system
was proved in \cite{EF70}, \cite[Chapter III]{FENO6} (see also
Bresch and Desjardins \cite{BRDE}, \cite{BRDE1} for an altrenative
approach based on a specific relation satisfied by the
density-dependent viscosity coefficients). Despite their apparent
success in the theory of \emph{existence}, the weak solutions in
most of the physically relevant cases are not (known to be)
uniquely determined by the data and may exhibit other rather
pathological properties, see Hoff and Serre \cite{HOSE}. 

A fundamental test of \emph{admissibility} of a class of weak
solutions to a given evolutionary problem is the property of
\emph{weak-strong uniqueness}. More specifically, the weak
solution must coincide with a (hypothetical) strong solution
emanating from the same initial data as long as the latter exists.
In other words, the strong solutions are unique within the class
of weak solutions. {This problem has been intensively studied for the \emph{incompressible} Navier-Stokes system 
for which, loosely speaking and simplifying, this is the ``best''
property that one can prove in this respect for the weak
solutions, see Prodi \cite{PR}, Serrin \cite{Serrin} for the first results in
this direction, and Escauriaza, Seregin, \v Sver\'ak \cite{ESSV}
for the most recent development.}

  The main goal of the present paper is
to show that the class of weak solutions to the \emph{full}
Navier-Stokes-Fourier system introduced in \cite{FENO6} enjoys
the weak-strong uniqueness property.

The motion of a general compressible, viscous, and heat conducting fluid is described by means of three basic state variables: the mass density $\vr=\vr(t,x)$, the velocity field
$\vu = \vu(t,x)$, and the absolute temperature $\vt = \vt(t,x)$, where $t$ is the time, and $x \in \Omega \subset R^3$ is the space variable in the Eulerian coordinate system. The time evolution of these quantities is governed by a system
of partial differential equations - mathematical formulation of the physical principles of balance of mass, momentum, and entropy:
\bFormula{i1}
\partial_t \vr + \Div (\vr \vu) = 0,
\eF
\bFormula{i2}
\partial_t (\vr \vu) + \Div (\vr \vu \otimes \vu) + \Grad p(\vr, \vt) = \Div \tn{S}(\vt, \Grad \vu),
\eF
\bFormula{i3}
\partial_t (\vr s(\vr, \vt)) + \Div (\vr s(\vr, \vt) \vu) + \Div \left( \frac{ \vc{q}(\vt, \Grad \vt)}{\vt} \right) = \sigma,
\eF
where $p$ is the pressure, $s$ is the (specific) entropy, and where we have deliberately ommited the influence of external sources. Furthermore, we suppose that the viscous stress  $\tn{S}$ is a linear function of the velocity gradient therefore described by Newton's law
\bFormula{i4}
\tn{S}(\vt, \Grad \vu) =  \mu(\vt) \Big( \Grad \vu + \Grad^t \vu - \frac{2}{3} \Div \vu \tn{I} \Big) + \eta (\vt) \Div \vu \tn{I},
\eF
while
$\vc{q}$ is the heat flux satisfying Fourier's law
\bFormula{i5}
\vc{q} = - \kappa(\vt) \Grad \vt,
\eF
and $\sigma$ stands for the entropy production rate specified below. The system of equations (\ref{i1} - \ref{i3}), with the constitutive relations (\ref{i4}), (\ref{i5}) is called \emph{Navier-Stokes-Fourier system}.
Equations (\ref{i1} - \ref{i4}) are supplemented with \emph{conservative} boundary conditions, say,
\bFormula{i6}
\vu|_{\partial \Omega} = \vc{q} \cdot \vc{n}|_{\partial \Omega} = 0.
\eF

A concept of \emph{weak solution} to the Navier-Stokes-Fourier
system (\ref{i1} - \ref{i6}) based on Second law of thermodynamics
was introduced in \cite{DUFE2}. The weak solutions satisfy the
field equations (\ref{i1} - \ref{i3}) in the sense of
distributions, where the entropy production $\sigma$ is a
non-negative measure, \bFormula{i7} \sigma \geq \frac{1}{\vt}
\left( \tn{S}(\vt, \Grad \vu) : \Grad \vu - \frac{\vc{q}(\vt,
\Grad \vt) \cdot \Grad \vt}{\vt} \right). \eF In order to
compensate for the lack of information resulting from the
inequality sign in (\ref{i7}), the resulting system is
supplemented by the total energy balance, \bFormula{i8} \frac{{\rm
d}}{{\rm d}t} \intO{ \left( \frac{1}{2} \vr |\vu|^2 + \vr e(\vr,
\vt) \right) } = 0, \eF where $e= e(\vr, \vt)$ is the (specific)
internal energy. Under these circumstances, it can be shown (see
\cite[Chapter 2]{FENO6}) that any weak solution of (\ref{i1}) that
is sufficiently \emph{smooth} satisfies, instead of (\ref{i7}),
the standard relation \bFormula{i9} \sigma = \frac{1}{\vt} \left(
\tn{S}(\vt, \Grad \vu) : \Grad \vu - \frac{\vc{q}(\vt, \Grad \vt)
\cdot \Grad \vt}{\vt} \right). \eF

Despite the fact that the framework of the weak solutions is sufficiently robust to provide:
\begin{itemize}
\item
a general \emph{existence theory} of global-in-time solutions under certain restrictions imposed on the
constitutive relations, see \cite[Chapter 3]{FENO6};
\item
quite satisfactory and complete description  of the long-time behavior of solutions, see \cite{FeiPr};
\item
rigorous justification of certain singular limits, see \cite{FENO6};
\end{itemize}
the class of weak solutions satisfying (\ref{i7}), (\ref{i8})
instead of (\ref{i9}) is probably too large to give rise to a
unique solution of the associated initial-boundary value problem.
The main issue seems to be the ambiguity of possible continuation
of solutions in the hypothetical presence of the so-called vacuum
zones (where $\vr = 0$) in combination with the possibility of
concentrations in the entropy production rate $\sigma$. On the
other hand, as we show in the present paper, any weak solution of
the Navier-Stokes-Fourier system coincides with a strong solution
as long as the latter exists. This is another piece of evidence
that the weak solutions, in the sense specified above, represent a
suitable extension of classical smooth solutions beyond their
(hypothetical) ``blow-up'' time.

Our approach is based on the method of \emph{relative entropy}
(cf. Carrillo et al. \cite{CaJuMaToUn}, Dafermos \cite{Daf4}, Saint-Raymond \cite{SaRay}), represented in the present context by the quantity
\bFormula{i10}
\mathcal{E}(\vr, \vt \  | \tilde \vr , \tilde \vt ) =
H_{\tilde \vt} (\vr, \vt) - \partial_\vr H_{\tilde \vt} (\tilde \vr, \tilde \vt) (\vr - \tilde \vr) - H_{\tilde \vt} (\tilde \vr, \tilde \vt),
\eF
where $H_{\tilde \vt}(\vr, \vt)$ is a thermodynamic potential termed \emph{ballistic free energy}
\bFormula{i11}
H_{\tilde \vt} (\vr, \vt) = \vr e(\vr, \vt) - \tilde \vt \vr s(\vr , \vt),
\eF
introduced by Gibbs and discussed more recently by Ericksen \cite{Eri}.

We assume that the thermodynamic functions $p$, $e$, and $s$ are interrelated through Gibbs' equation
\bFormula{i12}
\vt D s(\vr, \vt) = D e(\vr, \vt) + p(\vr, \vt) D \left( \frac{1}{\vr} \right).
\eF
The subsequent analysis leans essentially on \emph{thermodynamic stability} of the fluid system expressed through
\bFormula{i13}
\frac{\partial p(\vr, \vt)}{\partial \vr} > 0,\ \frac{\partial e(\vr, \vt)}{\partial \vt} > 0
\ \mbox{for all} \ \vr, \vt > 0.
\eF
In terms of the function $H_{\tilde \vt}$, thermodynamic stability implies that
\bFormula{i14}
\vr \mapsto H_{\tilde \vt} (\vr, \tilde \vt) \ \mbox{is strictly convex,}
\eF
while
\bFormula{i15}
\vt \mapsto H_{\tilde \vt}(\vr, \vt) \ \mbox{attains its global minimum at}\ \vt = \tilde \vt.
\eF
The above relations reflect \emph{stability} of the equilibrium solutions to the Navier-Stokes-Fourier system
(see Bechtel, Rooney, and Forest \cite{BEROFO}) and play a crucial role in the study of the long-time behavior
of solutions, see \cite[Chapters 5,6]{FeiPr}. Moreover, properties (\ref{i14}), (\ref{i5}) yield the necessary uniform bounds in
singular limits with ill-prepared initial data, see \cite[Chapter 2]{FENO6}. Last but not least, as we will see below,
the relations (\ref{i14}), (\ref{i15}) represent the key ingredient in the proof of weak-strong uniqueness.

The weak-strong uniqueness property for the standard
\emph{incompressible} Navier-Stokes system was established in
seminal papers by Prodi \cite{PR} and Serrin \cite{Serrin}. The
situation is a bit more delicate in the case of a
\emph{compressible} fluid. Germain \cite{Ger} proves weak-strong
uniqueness for the isentropic Navier-Stokes system, unfortunately,
in the class of ``more regular'' weak solutions which existence is
not known, see also Desjardin \cite{DES2} for the previous
results in this direction. The problem is finally settled in \cite{FeJiNo}
(see also \cite{FENOSU}), where the weak-strong uniqueness is
established in the class of the weak solutions to the barotropic
Navier-Stokes system satisfying the energy inequality.

Our method leans on a kind of  \emph{total dissipation balance} stated in terms of the quantity
\[
I = \intO{ \left( \frac{1}{2} \vr |\vu - \tilde \vu |^2 + \mathcal{E}(\vr, \vt | \tilde \vr, \tilde \vt ) \right)} ,
\]
where $\{ \vr, \vt, \vu \}$ is a weak solution to the Navier-Stokes-Fourier system and
$\{ \tilde \vr, \tilde \vt, \tilde \vu \}$ is a hypothetical classical solution emanating from the same initial data. We use a Gronwall type argument to show that $I(t) = 0$ provided
$I(0) = 0$ as long as the classical solution exists. The proof is carried over under the list of structural hypotheses
that are identical to those providing the existence of global-in-time solutions. Possible relaxations are discussed
in the concluding part of the paper.

The paper is organized as follows. In Section \ref{m}, we recall the definition of the weak solutions to the Navier-Stokes-Fourier system, introduce the hypotheses imposed on the
constitutive relations, and state our main result. In Section \ref{d}, we deduce a version of relative entropy inequality by taking $\tilde \vr = r$, $\tilde \vt = \Theta$, $\tilde \vu = \vc{U}$,
where $r$, $\Theta$, and $\vc{U}$ are arbitrary smooth functions satisfying the relevant boundary conditions. Finally, in Section \ref{w}, we establish the weak-strong uniqueness property. Possible extensions and further applications of the method are discussed in Section \ref{c}.

\section{Main result}
\label{m}

Motivated by the existence theory developed in \cite[Chapter 3]{FENO6}, we assume that the pressure $p = p(\vr, \vt)$ can be written in the form
\bFormula{m1}
p(\vr, \vt) = \vt^{5/2} P \left( \frac{\vr}{\vt^{3/2}} \right) + \frac{a}{3} \vt^4, \
a > 0,
\eF
where
\bFormula{m2}
P \in C^1[0, \infty), \ P(0) = 0, \ P'(Z) > 0 \ \mbox{for all}\ Z \geq 0.
\eF
In agreement with Gibbs' relation (\ref{i12}), the (specific) internal energy can be taken as
\bFormula{m3}
e(\vr, \vt) = \frac{3}{2} \frac{\vt^{5/2}}{\vr}P \left( \frac{\vr}{\vt^{3/2}} \right) +
a \frac{\vt^4}{\vr}.
\eF
Furthermore, by virtue of the second inequality in thermodynamic stability hypothesis (\ref{i13}), we have
\bFormula{m4}
0 < \frac{\frac{5}{3} P(Z) - P'(Z) Z }{Z} < c \ \mbox{for all}\ Z > 0.
\eF
Relation (\ref{m4}) implies that the function $Z \mapsto P(Z) / Z^{5/3}$ is decreasing, and we suppose that
\bFormula{m5}
\lim_{Z \to \infty} \frac{P(Z)}{Z^{5/3}} = P_\infty > 0.
\eF
Finally, the formula for (specific) entropy reads
\bFormula{m6}
s(\vr, \vt) = S \left( \frac{\vr}{\vt^{3/2}} \right)  + \frac{4a}{3} \frac{\vt^3}{\vr},
\eF
where, in accordance with Third law of thermodynamics,
\bFormula{m7}
S'(Z) = - \frac{3}{2} \frac{ \frac{5}{3} P(Z) - P'(Z) Z }{Z^2} < 0,\
\lim_{ Z \to \infty } S(Z) = 0.
\eF
The reader may consult Eliezer, Ghatak, and Hora \cite{EGH} and \cite[Chapter 3]{FENO6}
for the physical background and further discussion concerning the structural hypotheses
(\ref{m1} - \ref{m5}).

For the sake of simplicity and clarity of presentation, we take the transport coefficients in the form
\bFormula{m8}
\mu (\vt) = \mu_0 + \mu_1 \vt ,\ \mu_0, \mu_1 > 0, \ \eta \equiv 0,
\eF
\bFormula{m9}
\kappa (\vt) = \kappa_0 + \kappa_2 \vt^2 + \kappa_3 \vt^3 ,\ \kappa_i > 0, \ i=0,2,3,
\eF
cf. \cite[Chapter 3]{FENO6}.

\subsection{Weak solutions to the Navier-Stokes-Fourier system}
\label{ws}

Let $\Omega \subset R^3$ be a bounded Lipschitz domain. We say
that a trio $\{ \vr, \vt, \vu \}$ is a \emph{weak solution} to the
Navier-Stokes-Fourier system (\ref{i1} - \ref{i8}) emanating from
the initial data \bFormula{m10} \vr(0, \cdot) = \vr_0, \ \vr \vu
(0, \cdot) = \vr_0 \vu_0,\ \vr s(\vr, \vt)(0, \cdot) = \vr_0
s(\vr_0, \vt_0),\quad{ \vr_0\ge 0,\;\vt_0>0} \eF   if:

\begin{itemize}

\item the density and the absolute temperature satisfy $\vr(t,x) \geq 0$, $\vt(t,x) > 0$ for a.a. $(t,x) \in (0,T) \times \Omega$, $\vr
\in C_{\rm weak}([0,T]; L^{5/3})$, $\vr \vu \in C_{\rm weak}([0,T]; L^{5/4}(\Omega;R^3))$,
$\vt \in L^\infty(0,T; L^4(\Omega)) \cap L^2(0,T; W^{1,2}(\Omega))$, and $\vu \in L^2(0,T; W^{1,2}_0(\Omega;R^3))$;

\item equation (\ref{i1}) is replaced by a family of integral identities
\bFormula{m11} \intO{ \vr (\tau, \cdot) \varphi (\tau, \cdot)}
-\intO{ \vr_0 \varphi (0, \cdot) } = \int_0^\tau \intO{ \Big( \vr
\partial_t \varphi + \vr \vu \cdot \Grad \varphi \Big) } \ \dt \eF
for any $\varphi \in C^1([0,T] \times \Ov{\Omega})$, and any $\tau
\in [0,T]$;

\item momentum equation (\ref{i2}) is satisfied in the sense of distributions, specifically,
\bFormula{m12}
\intO{ \vr \vu (\tau, \cdot) \cdot \varphi (\tau, \cdot) } -
\intO{ \vr_0 \vu_0 \cdot \varphi (0, \cdot) }
\eF
\[
\int_0^\tau \intO{ \Big( \vr \vu \cdot \partial_t \varphi + \vr \vu \otimes \vu : \Grad \varphi + p(\vr, \vt) \Div \varphi - \tn{S}(\vt, \Grad \vu) : \Grad \varphi \Big) } \ \dt
\]
for any $\varphi \in C^1([0,T] \times \Ov{\Omega}; R^3)$, $\varphi|_{\partial \Omega} = 0$,
and any $\tau \in [0,T]$;
\item the entropy balance (\ref{i3}), (\ref{i7}) is replaced by a family of integral inequalities
\bFormula{m13}
\intO{ \vr_0 s(\vr_0, \vt_0) \varphi (0, \cdot) } -
\intO{ \vr s(\vr, \vt) (\tau, \cdot) \varphi (\tau, \cdot) }
\eF
\[
+
\int_0^\tau \intO{ \frac{ \varphi }{\vt} \left(
\tn{S}(\vt, \Grad \vu) : \Grad \vu - \frac{\vc{q}(\vt, \Grad \vt) \cdot \Grad \vt }{\vt}
\right) } \ \dt
\]
\[
\leq - \int_0^\tau \intO{ \left( \vr s(\vr, \vt) \partial_t \varphi + \vr s(\vr, \vt) \vu \cdot \Grad \varphi + \frac{ \vc{q}(\vt, \Grad \vt) \cdot \Grad \varphi }{\vt} \right)
} \ \dt
\]
for any $\varphi \in C^1([0,T] \times \Ov{\Omega})$, $\varphi \geq 0$, and a.a. $\tau \in [0,T]$;

\item the total energy is conserved:
\bFormula{m14}
\intO{ \left( \frac{1}{2} \vr |\vu|^2 + \vr e(\vr, \vt) \right)(\tau, \cdot) } =
\intO{ \left( \frac{1}{2} \vr_0 |\vu_0|^2 + \vr_0 e(\vr_0, \vt_0) \right) }
\eF
for a.a. $\tau \in [0,T]$.

\end{itemize}

The \emph{existence} of global-in-time weak solutions to the Navier-Stokes-Fourier system was established in \cite[Chapter 3, Theorem 3.1]{FENO6}.

\subsection{Main result}

We say that $\{ \tilde \vr, \tilde \vt, \tilde \vu \}$ is a
classical (strong) solution to the Navier-Stokes-Fourier system in
$(0,T) \times \Omega$ if
\bFormula{2.15} \tilde \vr \in C^1([0,T] \times \Ov{\Omega}),\
\tilde \vt,\ \partial_t \tilde \vt, \ \nabla^2 \tilde \vt \in
C([0,T] \times \Ov{\Omega}), \ \tilde \vu, \
\partial_t \tilde \vu,\ \nabla^2 \tilde \vu \in C([0,T] \times
\Ov{\Omega};R^3), \eF
\[
\tilde \vr (t,x) \geq \underline{\vr} > 0 , \ \tilde \vt (t,x) \geq \underline{\vt} > 0
\ \mbox{for all}\ (t,x),
\]
and $\tilde \vr$, $\tilde \vt$, $\tilde \vu$ satisfy equations
(\ref{i1} - \ref{i3}), (\ref{i9}), together with the boundary
conditions (\ref{i6}).
{Observe that hypothesis (\ref{2.15}) implies the following regularity properties of the initial data:
\bFormula{2.16} \vr(0)=\vr_0\in
C^1(\overline\Omega),\quad\vr_0\ge\underline\vr>0, \eF
\[
\vt(0)=\vt_0\in C^2(\overline\Omega),\quad\vt_0\ge\underline\vt>0,
\]
\[
\vu(0)=\vu_0\in C^2(\overline\Omega).
\]
}

We are ready to state the main result of this paper.

\bTheorem{m1} Let $\Omega \subset R^3$ be a bounded Lipschitz
domain. Suppose that the thermodynamic functions $p$, $e$, $s$
satisfy hypotheses (\ref{m1} - \ref{m7}), and that the transport
coefficients $\mu$, $\eta$, and $\kappa$ obey (\ref{m8}),
(\ref{m9}). Let $\{ \vr, \vt, \vu \}$ be a weak solution of the
Navier-Stokes-Fourier system in $(0,T) \times \Omega$ in the sense
specified in Section \ref{ws}, and let $\{ \tilde \vr, \tilde \vt,
\tilde \vu \}$ be a strong solution emanating from the same
initial data
{(\ref{2.16}).
}

Then
\[
\vr \equiv \tilde \vr, \ \vt \equiv \tilde \vt,\ \vu \equiv \tilde \vu.
\]

\eT

The rest of the paper is devoted to the proof of Theorem \ref{Tm1}. Possible generalizations are discussed in Section \ref{c}.

\section{Relative entropy inequality}
\label{d}

We deduce a relative entropy inequality satisfied by \emph{any} weak solution to the Navier-Stokes-Fourier system. To this end, consider a trio $\{ r , \Theta, \vc{U} \}$ of
smooth functions, $r$ and $\Theta$ bounded below away from zero in $[0,T] \times \Omega$, and
$\vc{U}|_{\partial \Omega} = 0$.

Taking $\varphi = \frac{1}{2} |\vc{U}|^2$ as a test function in (\ref{m11}), we get
\bFormula{d1}
\intO{ \frac{1}{2} \vr |\vc{U}|^2 (\tau, \cdot) } = \intO{ \frac{1}{2} \vr_0 |\vc{U}(0, \cdot)|^2 } +
\int_0^\tau \intO{ \Big( \vr \vc{U} \cdot \partial_t \vc{U} + \vr \vu \cdot \Grad \vc{U} \cdot \vc{U} \Big) } \ \dt.
\eF
Similarly, the choice $\varphi = \vc{U}$ in (\ref{m12}) gives rise to
\bFormula{d2}
\intO{ \vr \vu \cdot \vc{U} (\tau, \cdot) }
- \intO{ \vr_0 \vu_0 \cdot \vc{U}(0, \cdot) }
\eF
\[
= \int_0^\tau \intO{ \Big( \vr \vu \cdot \partial_t \vc{U} + \vr \vu \otimes \vu :  \Grad \vc{U} +
p(\vr, \vt) \Div \vc{U} - \tn{S}(\vt, \Grad \vu) : \Grad \vc{U} \Big) } \ \dt.
\]
Combining relations (\ref{d1}), (\ref{d2}) with the total energy balance (\ref{m14}) we may infer that
\bFormula{d3}
\intO{ \left( \frac{1}{2} \vr | \vu - \vc{U} |^2 + \vr e(\vr, \vt) \right)(\tau, \cdot) } =
\intO{ \left( \frac{1}{2} \vr_0 |\vu_0 - \vc{U}(0, \cdot) |^2 + \vr_0 e(\vr_0, \vt_0) \right) }
\eF
\[
+ \int_0^\tau \intO{ \left( \Big( \vr \partial_t \vc{U} + \vr \vu \cdot \Grad \vc{U} \Big) \cdot (\vc{U} - \vu) - p(\vr, \vt) \Div \vc{U} + \tn{S}(\vt, \Grad \vu) : \Grad \vc{U} \right) } \ \dt.
\]

Now, take $\varphi = \Theta > 0$ as a test function in the entropy inequality (\ref{m13})
to obtain
\bFormula{d4}
\intO{ \vr_0 s(\vr_0, \vt_0) \Theta(0, \cdot)  } -
\intO{ \vr s(\vr, \vt) \Theta (\tau, \cdot) }
\eF
\[
+
\int_0^\tau \intO{ \frac{\Theta}{\vt} \left(
\tn{S}(\vt, \Grad \vu): \Grad \vu - \frac{ \vc{q}(\vt, \Grad \vt) \cdot \Grad \vt }{\vt} \right) } \ \dt
\]
\[
\leq - \int_0^\tau \intO{ \left( \vr s(\vr, \vt) \partial_t \Theta + \vr s(\vr, \vt) \vu \cdot \Grad \Theta + \frac{ \vc{q}(\vt, \Grad \vt) }{\vt} \cdot \Grad \Theta \right) } \ \dt.
\]
Thus, the sum of (\ref{d3}), (\ref{d4}) reads
\bFormula{d5}
\intO{ \left( \frac{1}{2} \vr | \vu - \vc{U} |^2 + \vr e(\vr, \vt) -
\Theta \vr s(\vr, \vt)  \right)(\tau, \cdot) }
\eF
\[
+
\int_0^\tau \intO{ \frac{\Theta}{\vt} \left(
\tn{S}(\vt, \Grad \vu): \Grad \vu - \frac{ \vc{q}(\vt, \Grad \vt) \cdot \Grad \vt }{\vt} \right) } \ \dt
\]
\[
=
\intO{ \left( \frac{1}{2} \vr_0 |\vu_0 - \vc{U}(0, \cdot) |^2 + \vr_0 e(\vr_0, \vt_0)
- \Theta(0, \cdot) \vr_0 s(\vr_0, \vt_0)  \right) }
\]
\[
+ \int_0^\tau \intO{ \left( \Big( \vr \partial_t \vc{U} + \vr \vu \cdot \Grad \vc{U} \Big) \cdot (\vc{U} - \vu) - p(\vr, \vt) \Div \vc{U} + \tn{S}(\vt, \Grad \vu) : \Grad \vc{U} \right) } \ \dt
\]
\[
- \int_0^\tau \intO{ \left( \vr s(\vr, \vt) \partial_t \Theta + \vr s(\vr, \vt) \vu \cdot \Grad \Theta + \frac{ \vc{q}(\vt, \Grad \vt) }{\vt} \cdot \Grad \Theta \right) } \ \dt .
\]

Next, we take $\varphi = \partial_\vr H_\Theta(r, \Theta)$ as a test function in (\ref{m11}) to deduce that
\[
\intO{ \vr \partial_\vr H_\Theta (r, \Theta) (\tau, \cdot)} =
\intO{ \vr_0 \partial_\vr H_{\Theta(0, \cdot)} (r(0, \cdot), \Theta(0, \cdot) }
\]
\[
+ \int_0^\tau \intO{ \left(
\vr \partial_t \Big( \partial_\vr H_\Theta (r, \Theta) \Big) + \vr \vu \cdot
\Grad \Big( \partial_\vr H_\Theta (r, \Theta) \Big) \right) } \ \dt,
\]
which, combined with (\ref{d5}), gives rise to
\bFormula{d6}
\intO{ \left( \frac{1}{2} \vr | \vu - \vc{U}|^2 + H_\Theta (\vr, \vt) -
\partial_\vr (H_\Theta)(r, \Theta)(\vr - r) - H_\Theta (r, \Theta) \right)(\tau, \cdot)} +
\eF
\[
\int_0^\tau \intO{ \frac{\Theta}{\vt} \left(
\tn{S}(\vt, \Grad \vu): \Grad \vu - \frac{ \vc{q}(\vt, \Grad \vt) \cdot \Grad \vt }{\vt} \right) } \ \dt
\]
\[
\leq \intO{  \frac{1}{2} \vr_0 | \vu_0 - \vc{U}(0, \cdot)|^2 }
\]
\[
 + \intO{ \left( H_{\Theta(0, \cdot)} (\vr_0, \vt_0) -
\partial_\vr (H_{\Theta(0, \cdot)})(r(0, \cdot), \Theta(0, \cdot))(\vr_0 - r(0, \cdot)) - H_{\Theta(0, \cdot)} (r(0, \cdot), \Theta(0, \cdot)) \right)}
\]
\[
+ \int_0^\tau \intO{ \left( \Big( \vr \partial_t \vc{U} + \vr \vu \cdot \Grad \vc{U} \Big) \cdot (\vc{U} - \vu) - p(\vr, \vt) \Div \vc{U} + \tn{S}(\vt, \Grad \vu) : \Grad \vc{U} \right) } \ \dt
\]
\[
- \int_0^\tau \intO{ \left( \vr s(\vr, \vt) \partial_t \Theta + \vr s(\vr, \vt) \vu \cdot \Grad \Theta + \frac{ \vc{q}(\vt, \Grad \vt) }{\vt} \cdot \Grad \Theta \right) } \ \dt .
\]
\[
- \int_0^\tau \intO{ \left(
\vr \partial_t \Big( \partial_\vr H_\Theta (r, \Theta) \Big) + \vr \vu \cdot
\Grad \Big( \partial_\vr H_\Theta (r, \Theta) \Big) \right) } \ \dt
\]
\[
+ \int_0^\tau \intO{ \partial_t \Big( r \partial_\vr (H_\Theta) (r, \Theta) -
H_\Theta (r, \Theta) \Big) } \ \dt.
\]

Furthermore, seeing that
\[
\partial_y \left( \partial_\vr H_\Theta (r, \Theta) \right) = - s(r, \Theta) \partial_y \Theta - r \partial_\vr (r, \Theta) \partial_y \Theta + \partial^2_{\vr,\vr} H_\Theta (r, \Theta) \partial_y \vr + \partial^2_{\vr, \vt} H_\Theta (r, \Theta) \partial_y \Theta
\]
for $y=t,x$,
we may rewrite (\ref{d6}) in the form
\bFormula{d7}
\intO{ \left( \frac{1}{2} \vr | \vu - \vc{U}|^2 + H_\Theta (\vr, \vt) -
\partial_\vr (H_\Theta)(r, \Theta)(\vr - r) - H_\Theta (r, \Theta) \right)(\tau, \cdot)} +
\eF
\[
\int_0^\tau \intO{ \frac{\Theta}{\vt} \left(
\tn{S}(\vt, \Grad \vu): \Grad \vu - \frac{ \vc{q}(\vt, \Grad \vt) \cdot \Grad \vt }{\vt} \right) } \ \dt
\]
\[
\leq \intO{  \frac{1}{2} \vr_0 | \vu_0 - \vc{U}(0, \cdot)|^2 }
\]
\[
 + \intO{ \left( H_{\Theta(0, \cdot)} (\vr_0, \vt_0) -
\partial_\vr (H_{\Theta(0, \cdot)})(r(0, \cdot), \Theta(0, \cdot))(\vr_0 - r(0, \cdot)) - H_{\Theta(0, \cdot)} (r(0, \cdot), \Theta(0, \cdot)) \right)}
\]
\[
+ \int_0^\tau \intO{ \left( \Big( \vr \partial_t \vc{U} + \vr \vu \cdot \Grad \vc{U} \Big) \cdot (\vc{U} - \vu) - p(\vr, \vt) \Div \vc{U} + \tn{S}(\vt, \Grad \vu) : \Grad \vc{U} \right) } \ \dt
\]
\[
- \int_0^\tau \intO{ \left( \vr \Big( s (\vr, \vt) - s(r, \Theta) \Big) \partial_t \Theta
+ \vr \Big( s(\vr, \vt) - s(r, \Theta) \Big) \vu \cdot \Grad \Theta +
\frac{ \vc{q}(\vt, \Grad \vt) }{\vt} \cdot \Grad \Theta \right) } \ \dt .
\]
\[
+ \int_0^\tau \intO{ \vr \Big( r \partial_\vr s(r, \Theta) \partial_t \Theta +
r  \partial_\vr s(r, \Theta) \vu \cdot \Grad \Theta \Big) } \ \dt
\]
\[
- \int_0^\tau \intO{
\vr \Big( \partial^2_{\vr, \vr} (H_\Theta) (r, \Theta) \partial_t r +
\partial^2_{\vr, \vt} (H_\Theta) (r, \Theta) \partial_t \Theta \Big) } \ \dt
\]
\[
- \int_0^\tau \intO{ \vr \vu \cdot
\Big( \partial^2_{\vr, \vr}( H_\Theta) (r, \Theta) \Grad r +
\partial^2_{\vr, \vt} (H_\Theta) (r, \Theta) \Grad \Theta \Big) \Big) } \ \dt
\]
\[
+ \int_0^\tau \intO{ \partial_t \Big( r \partial_\vr (H_\Theta) (r, \Theta) -
H_\Theta (r, \Theta) \Big) } \ \dt.
\]

In order to simplify (\ref{d7}), we recall several useful identities that follow directly from Gibbs' relation (\ref{i12}):
\[
\partial^2_{\vr, \vr} (H_\Theta) (r, \Theta) = \frac{1}{r} \partial_\vr p (r, \Theta),
\]
\bFormula{sim}
r \partial_\vr s(r, \Theta) = - \frac{1}{r} \partial_\vt p (r, \Theta),
\eF
\[
\partial^2_{\vr, \vt} (H_\Theta) (r, \Theta) = \partial_\vr \left( \vr ( \vt - \Theta )
\partial_\vt s  \right) (r, \Theta) = (\vt - \Theta) \partial_\vr \Big( \vr \partial_\vt s(\vr, \vt) \Big) (r, \Theta) = 0,
\]
and
\[
r \partial_\vr (H_\Theta) (r, \Theta) - H_\Theta(r, \Theta) = p(r, \Theta).
\]

Thus relation (\ref{d7}) can be finally written in a more concise form
\bFormula{d8}
\intO{ \left( \frac{1}{2} \vr | \vu - \vc{U}|^2 + \mathcal{E}(\vr, \vt | r, \Theta) \right)(\tau, \cdot)} +
\int_0^\tau \intO{ \frac{\Theta}{\vt} \left(
\tn{S}(\vt, \Grad \vu): \Grad \vu - \frac{ \vc{q}(\vt, \Grad \vt) \cdot \Grad \vt }{\vt} \right) } \ \dt
\eF
\[
\leq \intO{ \left( \frac{1}{2} \vr_0 | \vu_0 - \vc{U}(0, \cdot)|^2 +
\mathcal{E} (\vr_0, \vt_0 | r(0, \cdot), \Theta (0, \cdot)) \right)}
\]
\[
\int_0^\tau \intO{ \vr (\vu - \vc{U}) \cdot \Grad \vc{U} \cdot (\vc{U} - \vu)} \ \dt
+ \int_0^\tau \intO{ \vr \Big( s(\vr, \vt) - s(r, \Theta) \Big) \Big( \vc{U} - \vu \Big)
\cdot \Grad \Theta } \ \dt
\]
\[
+ \int_0^\tau \intO{ \left( \vr \Big(  \partial_t \vc{U} +  \vc{U} \cdot \Grad \vc{U} \Big) \cdot (\vc{U} - \vu) - p(\vr, \vt) \Div \vc{U} + \tn{S}(\vt, \Grad \vu) : \Grad \vc{U} \right) } \ \dt
\]
\[
- \int_0^\tau \intO{ \left( \vr \Big( s (\vr, \vt) - s(r, \Theta) \Big) \partial_t \Theta
+ \vr \Big( s(\vr, \vt) - s(r, \Theta) \Big) \vc{U} \cdot \Grad \Theta +
\frac{ \vc{q}(\vt, \Grad \vt) }{\vt} \cdot \Grad \Theta \right) } \ \dt
\]
\[
+ \int_0^\tau \intO{ \left( \left( 1 - \frac{\vr}{r} \right) \partial_t p(r, \Theta) -
\frac{\vr}{r} \vu \cdot \Grad p(r, \Theta) \right) } \ \dt,
\]
where $\mathcal{E}$ was introduced in (\ref{i10}).

Formula (\ref{d8}) represents a kind of \emph{relative entropy inequality} in the spirit of
\cite{FeJiNo}, \cite{FENOSU}. Note that it is satisfied for \emph{any} trio of smooth functions
$\{r, \Theta, \vc{U}\}$ provided $\vc{U}$ vanishes on the boundary $\partial \Omega$, meaning
$\vc{U}$ is an admissible test function in the weak formulation of the momentum equation (\ref{m12}). Similar result can be obtained for other kinds of boundary conditions. The requirement on smoothness of $r$, $\Theta$, and $\vc{U}$ can be relaxed given the specific
integrability properties of the weak solution $\vr$, $\vt$, $\vu$. It is also worth-noting that (\ref{d8}) holds provided $p$, $e$, and $s$ satisfy only Gibbs' equation, the structural restrictions introduced in Section \ref{m} are not needed at this step.

\section{Weak-strong uniqueness}
\label{w}

In this section we finish the proof of Theorem \ref{Tm1} by applying the relative entropy inequality (\ref{d8}) to $r= \tilde \vr$, $\Theta = \tilde \vt$, and $\vc{U} = \tilde \vu$,
where $\{ \tilde \vr, \tilde \vt, \tilde \vu \}$ is a classical (smooth) solution of the Navier-Stokes-Fourier system such that
\[
\tilde \vr (0, \cdot) = \vr_0 , \ \tilde \vu (0, \cdot) = \vu_0,\ \tilde \vt (0, \cdot) =
\vt_0.
\]
Accordingly, the integrals depending on the initial values on the right-hand side of (\ref{d8}) vanish, and we apply a Gronwall type argument to deduce the desired result, namely,
\[
\vr \equiv \tilde \vr ,\ \vt \equiv \tilde \vt, \ \mbox{and}\ \vu \equiv \tilde \vu.
\]
Here, the hypothesis of thermodynamic stability formulated in (\ref{i13}) will play a crucial role.

\subsection{Preliminaries, notation}

Following \cite[Chapters 4,5]{FENO6} we introduce \emph{essential} and \emph{residual} component of each quantity appearing in (\ref{d8}). To begin, we choose positive constants
$\underline{\vr}$, $\Ov{\vr}$, $\underline{\vt}$, $\Ov{\vt}$ in such a way that
\[
0 < \underline{\vr} \leq \frac{1}{2} \min_{(t,x) \in [0,T] \times \Ov{\Omega}} \tilde \vr(t,x) \leq 2 \max_{(t,x) \in [0,T] \times \Ov{\Omega}} \tilde \vr (t,x) \leq
\Ov{\vr},
\]
\[
0 < \underline{\vt} \leq \frac{1}{2} \min_{(t,x) \in [0,T] \times \Ov{\Omega}} \tilde \vt(t,x) \leq 2 \max_{(t,x) \in [0,T] \times \Ov{\Omega}} \tilde \vt (t,x) \leq
\Ov{\vt}.
\]

In can be shown, as a consequence of the hypothesis of thermodynamic stability (\ref{i13}),
or, more specifically, of (\ref{i14}), (\ref{i15}), that
\bFormula{w1}
\mathcal{E} (\vr, \vt | \tilde \vr, \tilde \vt ) \geq c \left\{
\begin{array}{l} |\vr - \tilde \vr |^2 + |\vt - \tilde \vt |^2
\ \mbox{if} \ (\vr, \vt) \in [\underline{\vr}, \overline{\vr}]
\times [\underline{\vt}, \Ov{\vt}] \\ \\
1 + |\vr s(\vr, \vt) | + \vr e(\vr, \vt) \ \mbox{otherwise,}
\end{array}
\right.
\eF
whenever $[\tilde \vr, \tilde \vt] \in  [\underline{\vr}, \overline{\vr}]
\times [\underline{\vt}, \Ov{\vt}]$, where the constant $c$ depends only on
$\underline{\vr}$, $\Ov{\vr}$, $\underline{\vt}$, $\Ov{\vt}$ and the structural properties of the
thermodynamic functions $e$, $s$, see \cite[Chapter 3, Proposition 3.2]{FENO6}.

Now, each measurable function $h$ can be written as
\[
h = h_{\rm ess} + h_{\rm res},
\]
where
\[
h_{\rm ess}(t,x) = \left\{ \begin{array}{l} h(t,x) \ \mbox{if}\ (\vr(t,x), \vt(t,x)) \in  [\underline{\vr}, \overline{\vr}]
\times [\underline{\vt}, \Ov{\vt}] \\ \\ 0 \ \mbox{otherwise}
\end{array} \right. ,\ h_{\rm res} = h - h_{\rm ess}.
\]

\subsection{Relative entropy balance}

Taking $r = \tilde \vr$, $\Theta = \tilde \vt$, $\vc{U} = \tilde \vu$ in (\ref{d8}) and using the fact that the initial values coincide, we obtain
\bFormula{w2}
\intO{ \left( \frac{1}{2} \vr | \vu - \tilde \vu |^2 + \mathcal{E}(\vr, \vt | \tilde \vr, \tilde \vt ) \right)(\tau, \cdot)} +
\int_0^\tau \intO{ \frac{\tilde \vt}{\vt} \left(
\tn{S}(\vt, \Grad \vu): \Grad \vu - \frac{ \vc{q}(\vt, \Grad \vt) \cdot \Grad \vt }{\vt} \right) } \ \dt
\eF
\[
\leq \int_0^\tau \intO{ \vr |\vu - \tilde \vu|^2  |\Grad \tilde \vu| } \ \dt
+ \int_0^\tau \intO{ \vr \Big( s(\vr, \vt) - s(\tilde \vr, \tilde \vt) \Big) \Big( \tilde \vu - \vu \Big)
\cdot \Grad \tilde \vt } \ \dt
\]
\[
+ \int_0^\tau \intO{ \left( \vr \Big(  \partial_t \tilde \vu +  \tilde \vu \cdot \Grad \tilde \vu \Big) \cdot (\tilde \vu - \vu) - p(\vr, \vt) \Div \tilde \vu + \tn{S}(\vt, \Grad \vu) : \Grad \tilde \vu \right) } \ \dt
\]
\[
- \int_0^\tau \intO{ \left( \vr \Big( s (\vr, \vt) - s(\tilde \vr, \tilde \vt) \Big) \partial_t \tilde \vt
+ \vr \Big( s(\vr, \vt) - s(\tilde \vr, \tilde \vt) \Big) \tilde \vu \cdot \Grad \tilde \vt +  \frac{ \vc{q}(\vt, \Grad \vt) }{\vt} \cdot \Grad \tilde \vt \right) } \ \dt
\]
\[
+ \int_0^\tau \intO{ \left( \left( 1 - \frac{\vr}{\tilde \vr} \right) \partial_t p(\tilde \vr, \tilde \vt) -
\frac{\vr}{\tilde \vr} \vu \cdot \Grad p(\tilde \vr, \tilde \vt) \right) } \ \dt.
\]

In order to handle the integrals on the right-hand side of (\ref{w2}), we proceed by several steps:

\medskip

{\bf Step 1:}

We have
\bFormula{w3}
\intO{ \vr |\vu - \tilde \vu|^2  |\Grad \tilde \vu| } \ \dt \leq
2 \| \Grad \tilde \vu \|_{L^\infty(\Omega; R^{3 \times 3} )}
\intO{ \frac{1}{2} \vr | \vu - \tilde \vu |^2 }.
\eF

\medskip

{\bf Step 2:}
\[
\left| \intO{ \vr \Big( s(\vr, \vt) - s(\tilde \vr, \tilde \vt) \Big) \Big( \tilde \vu - \vu \Big)
\cdot \Grad \tilde \vt } \right|
\]
\[
\| \Grad \tilde \vt \|_{L^\infty(\Omega;R^3)} \left[
2 \Ov{\vr} \intO{ \left| \left[ s(\vr, \vt) - s(\tilde \vr, \tilde \vt) \right]_{\rm ess} \right|
|\vu - \tilde \vu | } + \intO{ \left| \left[ \vr \left( s(\vr, \vt) - s(\tilde \vr, \tilde \vt) \right) \right]_{\rm res} \right| |\vu - \tilde \vu| }
\right],
\]
where, by virtue of (\ref{w1}),
\bFormula{w4}
\intO{ \left| \left[ s(\vr, \vt) - s(\tilde \vr, \tilde \vt) \right]_{\rm ess} \right|
|\vu - \tilde \vu | } \leq \delta \| \vu - \tilde \vu \|^2_{L^2(\Omega; R^3)} +
c(\delta) \intO{ \mathcal{E}(\vr, \vt | \tilde \vr, \tilde \vt ) }
\eF
for any $\delta > 0$.

Similarly, by interpolation inequality,
\[
\intO{ \left| \left[ \vr \left( s(\vr, \vt) - s(\tilde \vr, \tilde \vt) \right) \right]_{\rm res} \right| |\vu - \tilde \vu| }
\]
\[
 \leq \delta \| \vu - \tilde \vu \|^2_{L^6(\Omega;R^3)} + c(\delta)
\left\| \left[ \vr \left( s(\vr, \vt) - s(\tilde \vr, \tilde \vt) \right) \right]_{\rm res} \right\|_{L^{6/5}(\Omega)}^2
\]
for any $\delta > 0$, where, furthermore,
\[
\left|
\left[ \vr \left( s(\vr, \vt) - s(\tilde \vr, \tilde \vt) \right) \right]_{\rm res}
\right| \leq c \left[ \vr + \vr s(\vr, \vt) \right]_{\rm res} .
\]
In accordance with (\ref{m6}), (\ref{m7}),
\bFormula{4.5--} 0 \leq \vr s(\vr, \vt) \leq c \left( \vt^3 + \vr
S \left( \frac{\vr}{\vt^{3/2}} \right) \right) , \eF
 where
 {
\bFormula{4.5-} \vr S \left( \frac{\vr}{\vt^{3/2}} \right) \leq
c(\vr + \vr [\log\vt]^+ + \vr |\log(\vr)| ). \eF
}
 On the other hand, as a
direct consequence of hypotheses (\ref{m3} - \ref{m5}), we get
\bFormula{w5} \vr e(\vr, \vt) \geq c( \vr^{5/3} + \vt^4 ). \eF
{Finally, we observe that (\ref{m13}-\ref{2.15}) imply
\bFormula{4.5+}
t \mapsto  \int_\Omega {\cal E}(\vr,\vt|\tilde\vr,\tilde\vt) \ \dx \in
 L^\infty(0,T).
 \eF

 Consequently, using (\ref{w1}), estimates (\ref{4.5--}-\ref{w5}) and the Holder inequality,
 we may
infer that
}
\bFormula{w6} \left\| \left[ \vr \left( s(\vr, \vt) -
s(\tilde \vr, \tilde \vt) \right) \right]_{\rm res}
\right\|_{L^{6/5}(\Omega)}^2 \leq c \left( \intO{ \mathcal{E}(\vr,
\vt | \tilde \vr, \tilde \vt ) } \right)^{5/3}. \eF

Combining (\ref{w4}), (\ref{w6})
{
with (\ref{4.5+}) we conclude that \bFormula{w7} \left| \intO{ \vr
\Big( s(\vr, \vt) - s(\tilde \vr, \tilde \vt) \Big) \Big( \tilde
\vu - \vu \Big) \cdot \Grad \tilde \vt } \right| \eF
\[
\leq \| \Grad \tilde \vt \|_{L^\infty(\Omega;R^3)} \left[ \delta
\| \vu - \tilde \vu \|^2_{W^{1,2}_0 (\Omega;R^3)} + c(\delta)
\intO{ \mathcal{E}(\vr, \vt | \tilde \vr, \tilde \vt ) }
\right]\le
\]
\[
\delta \| \vu - \tilde \vu \|^2_{W^{1,2}_0 (\Omega;R^3)} +
K(\delta,\cdot) \intO{ \mathcal{E}(\vr, \vt | \tilde \vr, \tilde
\vt ) }
\]
for any $\delta > 0$. {Here and hereafter, $K(\delta,\cdot)$
is a generic constant depending on $\delta$, $\tilde\vr$, $\tilde\vc u$,
$\tilde\vt$ through the norms induced by (\ref{2.15}),
(\ref{2.16}) and $\underline\vr$, $\underline\vt$, while $K(\cdot)$
is independent of $\delta$ but depends on $\tilde\vr$, $\tilde\vc
u$, $\tilde\vt$, $\underline\vr$, $\underline\vt$ through the
norms  (\ref{2.15}), (\ref{2.16}).}
}

\medskip

{\bf Step 3:}

Writing
\[
\intO{  \vr \Big(  \partial_t \tilde \vu +  \tilde \vu \cdot \Grad \tilde \vu \Big)
\cdot (\tilde \vu - \vu) } = \intO{ \frac{\vr}{\tilde \vr} (\tilde \vu - \vu) \cdot \left(
\Div \tn{S}(\tilde \vt, \Grad \tilde \vu) - \Grad p(\tilde \vr, \tilde \vt)  \right) }
\]
\[
= \intO{ \frac{1}{\tilde \vr} (\vr - \tilde \vr)  (\tilde \vu - \vu) \cdot \left(
\Div \tn{S}(\tilde \vt, \Grad \tilde \vu) - \Grad p(\tilde \vr, \tilde \vt)  \right) }
+ \intO{ (\tilde \vu - \vu) \cdot \left(
\Div \tn{S}(\tilde \vt, \Grad \tilde \vu) - \Grad p(\tilde \vr, \tilde \vt)  \right) },
\]
we observe that the first integral on the right-hand side can be
handled in the same way as in Step 2,
{
namely,
$$
 \intO{ \left[\frac{1}{\tilde \vr} (\vr - \tilde \vr)  (\tilde \vu - \vu) \cdot \left(
\Div \tn{S}(\tilde \vt, \Grad \tilde \vu) - \Grad p(\tilde \vr,
\tilde \vt)  \right)\right]_{\rm ess} }\le
K(\delta, \cdot)\left\|\left[\vr-\tilde\vr\right]_{\rm
ess}\right\|^2_{L^2(\Omega)}+\delta\|\vc u-\tilde\vc
u\|^2_{L^2(\Omega;R^3)},
$$
$$
 \intO{ \left[\frac{1}{\tilde \vr} (\vr - \tilde \vr)  (\tilde \vu - \vu) \cdot \left(
\Div \tn{S}(\tilde \vt, \Grad \tilde \vu) - \Grad p(\tilde \vr,
\tilde \vt)  \right)\right]_{\rm res }}
$$
$$
\le
K(\delta, \cdot)\left(\left\|\left[\vr\right]_{\rm
res}\right\|_{L^{6/5}(\Omega)}^2 + \left\|\left[1\right]_{\rm
res}\right\|^2_{L^{6/5}(\Omega)}\right)+\delta\|\vc u-\tilde\vc
u\|^2_{L^6(\Omega;R^3)},
$$
while, integrating by parts,
\[
\intO{ (\tilde \vu - \vu) \cdot \left(
\Div \tn{S}(\tilde \vt, \Grad \tilde \vu) - \Grad p(\tilde \vr, \tilde \vt)  \right) }
\]
\[
= \intO{ \left(
\tn{S}(\tilde \vt, \Grad \tilde \vu) : \Grad (\vu - \tilde \vu)  + p(\tilde \vr, \tilde \vt)
\Div (\tilde \vu - \vu )   \right) }
\]
Thus using again (\ref{w1}), (\ref{4.5+}) and continuous imbedding
$W^{1,2}(\Omega)\hookrightarrow L^6(\Omega)$ we arrive at
}
\bFormula{w8} \left| \intO{  \vr \Big(
\partial_t \tilde \vu +  \tilde \vu \cdot \Grad \tilde \vu \Big)
\cdot (\tilde \vu - \vu) } \right| \leq \intO{ \left(
\tn{S}(\tilde \vt, \Grad \tilde \vu) : \Grad (\vu - \tilde \vu)  +
p(\tilde \vr, \tilde \vt) \Div (\tilde \vu - \vu )   \right) } \eF
\[
+ c\left( |\Grad \tilde \vr|, |\Grad \tilde \vt |, |\nabla^2_x
\tilde \vu | \right) \left[ \delta \| \vu - \tilde \vu
\|^2_{W^{1,2}_0(\Omega;R^3)} + c(\delta) \intO{ \mathcal{E}(\vr,
\vt | \tilde \vr, \tilde \vt ) }   \right]\le
\]
{
\[
\intO{ \left( \tn{S}(\tilde \vt, \Grad \tilde \vu) : \Grad (\vu -
\tilde \vu)  + p(\tilde \vr, \tilde \vt) \Div (\tilde \vu - \vu )
\right) }
\]
\[
+ \delta \| \vu - \tilde \vu \|^2_{W^{1,2}(\Omega;R^3)}
+ K(\delta,\cdot) \intO{ \mathcal{E}(\vr, \vt | \tilde \vr, \tilde
\vt ) }
\]
}
 for any $\delta > 0$.

\medskip

{\bf Step 4:}

Next, we get
\[
\intO{ \vr \Big( s(\vr, \vt ) - s(\tilde \vr, \tilde \vt) \Big) \partial_t \tilde \vt } =
\]
\[
\intO{ \vr \Big[ s(\vr, \vt ) - s(\tilde \vr, \tilde \vt) \Big]_{\rm ess} \partial_t \tilde \vt }  + \intO{ \vr \Big[ s(\vr, \vt ) - s(\tilde \vr, \tilde \vt) \Big]_{\rm res} \partial_t \tilde \vt },
\]
where
{
\bFormula{w9} \left| \intO{ \vr \Big[ s(\vr, \vt ) -
s(\tilde \vr, \tilde \vt) \Big]_{\rm res} \partial_t \tilde \vt }
\right|
\eF
\[
\leq \| \partial_t \tilde \vt
\|_{L^\infty(\Omega)}
\Big(\int_\Omega\left[\vr s(\vr,\vt)\right]_{\rm res}{\rm d}x
+\|s(\tilde\vr,\tilde\vt)\|_{L^\infty(\Omega)}\int_\Omega[\vr]_{\rm
res}{\rm d}x\Big)\le
 K(\cdot) \intO{ \mathcal{E}(\vr, \vt| \tilde
\vr, \tilde \vt ) }, \]
}
while
\[
\intO{ \vr \Big[ s(\vr, \vt ) - s(\tilde \vr, \tilde \vt) \Big]_{\rm ess} \partial_t \tilde \vt }
\]
\[
 = \intO{ ( \vr - \tilde \vr)  \Big[ s(\vr, \vt ) - s(\tilde \vr, \tilde \vt) \Big]_{\rm ess} \partial_t \tilde \vt } + \intO{ \tilde \vr \Big[ s(\vr, \vt ) - s(\tilde \vr, \tilde \vt) \Big]_{\rm ess} \partial_t \tilde \vt },
\]
where,
{
with help of Taylor-Lagrange formula,
\bFormula{w10} \left| \intO{ ( \vr - \tilde \vr) \Big[ s(\vr, \vt
) - s(\tilde \vr, \tilde \vt) \Big]_{\rm ess}
\partial_t \tilde \vt } \right|
\eF
\[
\leq
 \Big(\sup_{(\vr,\vt)\in
[\underline\vr,\overline\vr]\times[\underline\vt,\overline\vt]}
|\partial_\vr s(\vr,\vt)|+ \sup_{(\vr,\vt)\in
[\underline\vr,\overline\vr]\times[\underline\vt,\overline\vt]}|\partial_\vt
s(\vr,\vt)|\Big)\;\|
\partial_t \tilde \vt
\|_{L^\infty(\Omega)}\times \]
\[
\times \int_\Omega\Big[\Big|[\vr-\tilde\vr]_{\rm
ess}\Big|\Big(\Big|[\vr-\tilde\vr]_{\rm ess}\Big|+
\Big|[\vt-\tilde\vt]_{\rm ess}\Big| \Big)\Big]{\rm d}x \le
K(\cdot)\intO{ \mathcal{E}(\vr, \vt| \tilde \vr, \tilde \vt ) }.
\]
}

Finally, we write
\[
\intO{ \tilde \vr \Big[ s(\vr, \vt ) - s(\tilde \vr, \tilde \vt) \Big]_{\rm ess} \partial_t \tilde \vt } = \intO{ \tilde \vr \Big[ s(\vr, \vt) - \partial_\vr s(\tilde \vr, \tilde \vt)
(\vr - \tilde \vr) - \partial_\vt s(\tilde \vr, \tilde \vt)
(\vt - \tilde \vt) - s( \tilde \vr, \tilde \vt) \Big]_{\rm ess} \partial_t \tilde \vt }
\]
\[
- \intO{ \tilde \vr \Big[ \partial_\vr s(\tilde \vr, \tilde \vt)
(\vr - \tilde \vr) + \partial_\vt s(\tilde \vr, \tilde \vt)
(\vt - \tilde \vt) \Big]_{\rm res} \partial_t \tilde \vt } + \intO{ \tilde \vr \Big[ \partial_\vr s(\tilde \vr, \tilde \vt)
(\vr - \tilde \vr) + \partial_\vt s(\tilde \vr, \tilde \vt)
(\vt - \tilde \vt) \Big] \partial_t \tilde \vt },
\]
where the first two integrals on the right-hand side can be
estimated exactly as in (\ref{w9}), (\ref{w10}).
{Thus we conclude that \bFormula{w11} - \intO{ \vr \Big( s(\vr, \vt
) - s(\tilde \vr, \tilde \vt) \Big) \partial_t \tilde \vt } \leq
K(\cdot) \intO{ \mathcal{E}(\vr, \vt| \tilde \vr, \tilde \vt ) }
\eF
}
\[
- \intO{ \tilde \vr \Big[ \partial_\vr s(\tilde \vr, \tilde \vt)
(\vr - \tilde \vr) + \partial_\vt s(\tilde \vr, \tilde \vt)
(\vt - \tilde \vt) \Big] \partial_t \tilde \vt }.
\]

\medskip

{\bf Step 5:}

Similarly to Step 4, we get \bFormula{w12} - \intO{ \vr \Big(
s(\vr, \vt ) - s(\tilde \vr, \tilde \vt) \Big) \tilde \vu \cdot
\Grad \tilde \vt } \leq K(\cdot) \intO{ \mathcal{E}(\vr, \vt|
\tilde \vr, \tilde \vt ) } \eF
\[
- \intO{ \tilde \vr \Big[ \partial_\vr s(\tilde \vr, \tilde \vt)
(\vr - \tilde \vr) + \partial_\vt s(\tilde \vr, \tilde \vt) (\vt -
\tilde \vt) \Big] \tilde \vu \cdot \Grad \tilde \vt }.
\]

\medskip

{\bf Step 6:}

Finally, we have
\bFormula{w13}
\intO{ \left( \left( 1 - \frac{\vr}{\tilde \vr} \right) \partial_t p(\tilde \vr, \tilde \vt)
- \frac{\vr}{\tilde \vr} \vu \cdot \Grad p(\tilde \vr, \tilde \vt) \right) }
\eF
\[
= \intO{ \left(\tilde \vr - \vr \right) \frac{1}{\tilde \vr} \left(
\partial_t p(\tilde \vr, \tilde \vt) +  \tilde \vu \cdot \Grad p(\tilde \vr, \tilde \vt) \right) } + \intO{ p(\tilde \vr, \tilde \vt) \Div \vu }
\]
\[
+ \intO{ (\vr - \tilde \vr) \frac{1}{\tilde \vr} \Grad p(\tilde \vr, \tilde \vt)\cdot
(\vu - \tilde \vu) },
\]
where, by means of the same arguments as in Step 2,
\[ \left|
\intO{ (\vr - \tilde \vr) \frac{1}{\tilde \vr} \Grad p(\tilde \vr,
\tilde \vt)\cdot (\vu - \tilde \vu) } \right|\le \]
\[
 c\left( |\Grad \tilde \vr|, |\Grad \tilde \vt| \right) \left[\delta \| \vu - \tilde \vu \|^2_{W^{1,2}(\Omega;R^3)} +
 \intO{ \mathcal{E}(\vr, \vt | \tilde \vr, \tilde \vt ) } \right]
\]
for any $\delta > 0$.
{
Resuming this step, we have \bFormula{w14} \intO{ \left( \left( 1
- \frac{\vr}{\tilde \vr} \right) \partial_t p(\tilde \vr, \tilde
\vt) - \frac{\vr}{\tilde \vr} \vu \cdot \Grad p(\tilde \vr, \tilde
\vt) \right) } \le \eF
\[
 \intO{ \left(\tilde \vr - \vr \right) \frac{1}{\tilde \vr}
\left(
\partial_t p(\tilde \vr, \tilde \vt) +  \tilde \vu \cdot \Grad p(\tilde \vr, \tilde \vt) \right) }
+ \intO{ p(\tilde \vr, \tilde \vt) \Div \vu } +
\]
\[
 \delta \| \vu - \tilde \vu \|^2_{W^{1,2}(\Omega;R^3)} +
K(\delta,\cdot) \intO{ \mathcal{E}(\vr, \vt | \tilde \vr, \tilde
\vt ) } .
\]
}

\medskip

{\bf Step 7:}

Summing up the estimates (\ref{w3}), (\ref{w7}), (\ref{w8}), (\ref{w11} - \ref{w14}),
we can rewrite the relative entropy inequality (\ref{w2}) in the form
\bFormula{w15}
\intO{ \left( \frac{1}{2} \vr |\vu - \tilde \vu |^2
+ \mathcal{E}(\vr, \vt| \tilde \vr, \tilde \vt ) \right) (\tau, \cdot) }
\eF
\[
+ \int_0^\tau \intO{ \left( \frac{\tilde \vt}{\vt} \tn{S} (\vt, \Grad \vu) : \Grad \vu -
\tn{S} (\tilde \vt, \Grad \tilde \vu ): (\Grad \vu - \Grad \tilde \vu )  -
\tn{S}(\vt, \Grad \vu) : \Grad \tilde \vu  \right) } \ \dt
\]
\[
+ \int_0^\tau \intO{ \left( \frac{\vc{q} (\vt, \Grad \vt) \cdot
\Grad \tilde \vt }{\vt} - \frac{\tilde \vt}{\vt}  \frac{\vc{q}
(\vt, \Grad \vt) \cdot \Grad  \vt }{\vt} \right) } \ \dt\leq
\]
{
\[
\int_0^\tau \left[ \delta \| \vu - \tilde \vu
\|^2_{W^{1,2}(\Omega;R^3)} + K(\delta,\cdot)\intO{ \left(
\frac{1}{2} \vr | \vu - \tilde \vu |^2 + \mathcal{E}(\vr, \vt |
\tilde \vr, \tilde \vt ) \right) }  \right] \ \dt
\]
}
\[
+ \int_0^\tau \intO{ \Big( p(\tilde \vr, \tilde \vt) - p(\vr, \vt) \Big) \Div \tilde \vu }
\ \dt
\]
\[
+ \int_0^\tau \intO{ (\tilde \vr - \vr) \frac{1}{\tilde \vr} \left[ \partial_t p(\tilde \vr, \tilde \vt) + \tilde \vu \cdot \Grad p(\tilde \vr, \tilde \vt) \right] } \ \dt
\]
\[
- \int_0^\tau \intO{ \tilde \vr \left( \partial_\vr s(\tilde \vr, \tilde \vt) (\vr - \tilde \vr) + \partial_\vt s(\tilde \vr, \tilde \vt) (\vt - \tilde \vt) \right) \left[
\partial_t \tilde \vt + \tilde \vu \cdot \Grad \tilde \vt \right] } \ \dt
\]
for any $\delta > 0$.

\medskip

{\bf Step 8:}

Our next goal is to control the last three integrals on the right-hand side of (\ref{w15}).
To this end, we use (\ref{sim}) to obtain
\[
\intO{ (\tilde \vr - \vr) \frac{1}{\tilde \vr} \left[ \partial_t p(\tilde \vr, \tilde \vt) + \tilde \vu \cdot \Grad p(\tilde \vr, \tilde \vt) \right] }
\]
\[
- \intO{ \tilde \vr \left( \partial_\vr s(\tilde \vr, \tilde \vt) (\vr - \tilde \vr) + \partial_\vt s(\tilde \vr, \tilde \vt) (\vt - \tilde \vt) \right) \left[
\partial_t \tilde \vt + \tilde \vu \cdot \Grad \tilde \vt \right] }
\]
\[
= \intO{ \tilde \vr (\tilde \vt - \vt) \partial_\vt s(\tilde \vr, \tilde \vt)
\left[ \partial_t \tilde \vt + \tilde \vu \cdot \Grad \tilde \vt \right] } +
\intO{ (\tilde \vr - \vr ) \frac{1}{\tilde \vr} \partial_\vr p(\tilde \vr, \tilde \vt)
\left[ \partial_t \tilde \vr + \tilde \vu \cdot \Grad \tilde \vr \right] },
\]
where, as $\tilde \vr$, $\tilde \vu$ satisfy the equation of continuity (\ref{i1}),
\bFormula{w16}
\intO{ (\tilde \vr - \vr ) \frac{1}{\tilde \vr} \partial_\vr p(\tilde \vr, \tilde \vt)
\left[ \partial_t \tilde \vr + \tilde \vu \cdot \Grad \tilde \vr \right] } = -
\intO{ (\tilde \vr - \vr) \partial_\vr p(\tilde \vr, \tilde \vt) \Div \tilde \vu }.
\eF

Finally, using (\ref{sim}) once more, we deduce that
\bFormula{w17}
\intO{ \tilde \vr (\tilde \vt - \vt) \partial_\vt s(\tilde \vr, \tilde \vt)
\left[ \partial_t \tilde \vt + \tilde \vu \cdot \Grad \tilde \vt \right] }
\eF
\[
\intO{ \tilde \vr (\tilde \vt - \vt)
\left[ \partial_t s(\tilde \vr, \tilde \vt) + \tilde \vu \cdot \Grad s(\tilde \vr, \tilde \vt) \right] } - \intO{ (\tilde \vt - \vt)
\partial_\vt p(\tilde \vr, \tilde \vt) \Div \tilde \vu }
\]
\[
= \intO{ (\tilde \vt - \vt) \left[ \frac{1}{\tilde \vt} \left(
\tn{S} (\tilde \vt, \Grad \tilde \vu) : \Grad \tilde \vu - \frac{\vc{q}(\tilde \vt, \Grad \tilde \vt) \cdot \Grad \tilde \vt}{\tilde \vt} \right) - \Div
\left( \frac{\vc{q} (\tilde \vt, \Grad \tilde \vt)}{\tilde \vt} \right) \right] }
\]
\[
- \intO{ (\tilde \vt - \vt)
\partial_\vt p(\tilde \vr, \tilde \vt) \Div \tilde \vu }
\]

{Seeing that
\[
\left| \intO{ \left( p(\tilde \vr, \tilde \vt) - \partial_\vr
p(\tilde \vr, \tilde \vt) (\tilde \vr - \vr) - \partial_\vt
p(\tilde \vr, \tilde \vt) (\tilde \vt - \vt) - p(\vr, \vt) \right)
\Div \tilde\vu } \right|
\]
}
\[
\leq
c \| \Div \tilde \vu \|_{L^\infty (\Omega)} \intO{ \mathcal{E} (\vr, \vt | \tilde \vr, \tilde \vt ) }
\]
we may use the previous relations to rewrite (\ref{w15}) in the form
\bFormula{w18}
\intO{ \left( \frac{1}{2} \vr |\vu - \tilde \vu |^2 + \mathcal{E}(\vr, \vt| \tilde \vr, \tilde \vt ) \right) (\tau, \cdot) }
\eF
\[
+ \int_0^\tau \intO{ \left( \frac{\tilde \vt}{\vt} \tn{S} (\vt, \Grad \vu) : \Grad \vu -
\tn{S} (\tilde \vt, \Grad \tilde \vu ): (\Grad \vu - \Grad \tilde \vu )  -
\tn{S}(\vt, \Grad \vu) : \Grad \tilde \vu  - \frac{\tilde \vt - \vt}{\tilde \vt}
\tn{S}(\tilde \vt, \Grad \tilde \vu ) : \Grad \tilde \vu \right) } \ \dt
\]
\[
+ \int_0^\tau \intO{ \left( \frac{\vc{q} (\vt, \Grad \vt) \cdot \Grad \tilde \vt }{\vt} -
\frac{\tilde \vt}{\vt}  \frac{\vc{q} (\vt, \Grad \vt) \cdot \Grad  \vt }{\vt}
+ (\tilde \vt - \vt) \frac{ \vc{q}(\tilde \vt, \Grad \tilde \vt) \cdot \Grad \tilde \vt }{
{\tilde \vt}^2 } + \frac{\vc{q}(\tilde \vt, \Grad \tilde \vt)}{\tilde \vt} \cdot
\Grad (\vt - \tilde \vt)
\right) }
\ \dt
\]
{
\[
\leq \int_0^\tau \left[ \delta \| \vu - \tilde \vu
\|^2_{W^{1,2}(\Omega;R^3)} + K(\delta,\cdot) \intO{ \left(
\frac{1}{2} \vr | \vu - \tilde \vu |^2 + \mathcal{E}(\vr, \vt |
\tilde \vr, \tilde \vt ) \right) }   \right] \ \dt
\]
}
 for any $\delta > 0$.

\subsection{Dissipative terms}

Our ultimate goal in the proof of Theorem \ref{Tm1} is to show that the ``dissipative'' terms appearing on the left-hand side of (\ref{w18}) containing $\Grad \vu$, $\Grad \vt$ are strong enough to control the $W^{1,2}-$norm of the velocity.

\subsubsection{Viscosity}

In accordance with hypothesis (\ref{m8}), we have
\[
\tn{S}(\vt, \Grad \vu) = \tn{S}^0 (\vt, \Grad \vu) + \tn{S}^1(\vt, \Grad \vu),
\]
where
\[
\tn{S}^0 (\vt, \Grad \vu) = \mu_0 \Big(\Grad \vu + \Grad^t \vu - \frac{2}{3} \Div \vu \tn{I} \Big), \ \tn{S}^1 (\vt, \Grad \vu) = \mu_1 \vt \Big(\Grad \vu + \Grad^t \vu - \frac{2}{3} \Div \vu \tn{I} \Big).
\]

Now, we write
\[
\frac{\tilde \vt}{\vt} \tn{S}^1 (\vt, \Grad \vu) : \Grad \vu -
\tn{S}^1 (\tilde \vt, \Grad \tilde \vu ): (\Grad \vu - \Grad \tilde \vu )  -
\tn{S}^1(\vt, \Grad \vu) : \Grad \tilde \vu  - \frac{\tilde \vt - \vt}{\tilde \vt}
\tn{S}^1(\tilde \vt, \Grad \tilde \vu ) : \Grad \tilde \vu
\]
\[
= \tilde \vt \left( \frac{\tn{S}^1 (\vt, \Grad \vu) }{\vt} - \frac{\tn{S}^1 (\tilde \vt,
\Grad \tilde \vu) }{\tilde \vt} \right) : (\Grad \vu - \Grad \tilde \vu )
+ (\tilde \vt - \vt ) \left( \frac{\tn{S}^1 (\vt, \Grad \vu) }{\vt} - \frac{\tn{S}^1 (\tilde \vt,
\Grad \tilde \vu) }{\tilde \vt} \right) : \Grad \tilde \vu,
\]
where, by virtue of Korn's inequality, \bFormula{w19} \intO{
\tilde \vt \left( \frac{\tn{S}^1 (\vt, \Grad \vu) }{\vt} -
\frac{\tn{S}^1 (\tilde \vt, \Grad \tilde \vu) }{\tilde \vt}
\right) : (\Grad \vu - \Grad \tilde \vu ) } \geq c \mu_1 \| \vu -
\tilde \vu \|^2_{W^{1,2}(\Omega;R^3)}. \eF

On the other hand, similarly to the preceding part, we can show that
\bFormula{w20}
\left| \intO{ (\tilde \vt - \vt ) \left( \frac{\tn{S}^1 (\vt, \Grad \vu) }{\vt} - \frac{\tn{S}^1 (\tilde \vt,
\Grad \tilde \vu) }{\tilde \vt} \right) : \Grad \tilde \vu } \right|
\eF
\[
\leq \| \Grad \tilde \vu \|_{L^\infty(\Omega; R^{3 \times 3})}
\left( \delta \| \vu - \tilde \vu \|^2_{W^{1,2} (\Omega;R^3)} +
c(\delta) \Big(\left\| [\vt-\tilde\vt]_{\rm
ess}\right\|^2_{L^2(\Omega)}+ \left\| [\vt-\tilde\vt]_{\rm
res}\right\|^2_{L^2(\Omega)}\Big)\right)
\]
\[
\leq \delta \| \vu - \tilde \vu \|^2_{W^{1,2} (\Omega;R^3)} +
K(\delta,\cdot) \intO{ \mathcal{E}(\vr, \vt | \tilde \vr, \tilde
\vt ) }
\]
for any $\delta > 0$, where we have used again (\ref{w1}).

{Next,
\[
\frac{\tilde \vt}{\vt} \tn{S}^0 (\Grad \vu) : \Grad \vu -
\tn{S}^0 (\Grad \tilde \vu ): (\Grad \vu - \Grad \tilde \vu )  -
\tn{S}^0(\Grad \vu) : \Grad \tilde \vu  - \frac{\tilde \vt - \vt}{\tilde \vt}
\tn{S}^0(\Grad \tilde \vu ) : \Grad \tilde \vu
\]
\[
= \frac{\tilde \vt}{\vt} \Big( \tn S^0 (\Grad \vu) - \tn S^0
(\Grad \tilde \vu) \Big) : \Grad (\vu - \tilde \vu ) + \tilde \vt
\left( \frac{1}{\vt} - \frac{1}{\tilde \vt} \right) \tn{S}^0
(\Grad \tilde \vu) : \Grad (\vu - \tilde \vu).
\]
%
%
First suppose that  $\vt \geq \tilde \vt$. Since the function $\vt
\mapsto 1/\vt$ is Lipschitz on the set $\vt \geq \tilde \vt$, we
conclude that
\bFormula{w21} \int_{ \{ \vt \geq \tilde \vt \} } \left| \tilde
\vt \left( \frac{1}{\vt} - \frac{1}{\tilde \vt} \right) \tn{S}^0
(\Grad \tilde \vu) : \Grad (\vu - \tilde \vu) \right| + \left|
\frac{ \tilde \vt - \vt }{\vt} \Big( \tn S^0(\Grad \vu) - \tn
S^0(\Grad \tilde \vu) \Big) : \Grad \tilde \vu \right| \ \dx \eF
\[
\leq \Big(\frac 1{\underline
\vt}+\|\Grad\tilde\vu\|_{L^\infty(\Omega;\R^3)}+\|\tilde\vt\|_{L^\infty(\Omega)}\Big)
\Big( \delta \| \Grad \vu - \Grad \tilde \vu \|^2_{L^2(\Omega;R^{3
\times 3})}+c(\delta)\|\vt-\tilde\vt\|^2_{L^2(\Omega)}\Big)
\]
\[
\delta \| \Grad \vu - \Grad \tilde \vu \|^2_{L^2(\Omega;R^{3
\times 3})} + K(\delta,\cdot) \intO{ \mathcal{E}(\vr, \vt | \tilde
\vr, \tilde \vt ) }
\]
for any $\delta > 0$.
}

 Finally, if $0 < \vt \leq \tilde \vt$, we have
\[
\frac{\tilde \vt}{\vt} \tn{S}^0 (\Grad \vu) : \Grad \vu -
\tn{S}^0 (\Grad \tilde \vu ): (\Grad \vu - \Grad \tilde \vu )  -
\tn{S}^0(\Grad \vu) : \Grad \tilde \vu  - \frac{\tilde \vt - \vt}{\tilde \vt}
\tn{S}^0(\Grad \tilde \vu ) : \Grad \tilde \vu
\]
\[
\geq \left( \tn{S}^0 (\Grad \vu) - \tn{S}^0 (\Grad \tilde \vu) \right): (\Grad \vu - \Grad \tilde \vu) + \frac{\tilde \vt - \vt}{\tilde \vt} \left[ \tn{S}^0 (\Grad \vu) : \Grad \vu -
\tn{S}^0 (\Grad \tilde \vu) : \Grad \tilde \vu \right];
\]
whence, by means of convexity of the function $\Grad \vu \mapsto \tn{S}^0 (\Grad \vu) : \Grad \vu$,
\[
\frac{\tilde \vt - \vt}{\tilde \vt} \left[ \tn{S}^0 (\Grad \vu) :
\Grad \vu - \tn{S}^0 (\Grad \tilde \vu) : \Grad \tilde \vu \right]
\geq \frac{\tilde \vt - \vt}{\tilde \vt} \tn{S}^0 (\Grad \tilde
\vu) \cdot \Grad (\vu - \tilde \vu),
\]
{
where, similarly as in (\ref{w21}),
\[
\int_\Omega\frac{\tilde \vt - \vt}{\tilde \vt} \tn{S}^0 (\Grad
\tilde \vu) \cdot \Grad (\vu - \tilde \vu){\rm d}x \le
K(\delta,\cdot) \intO{ \mathcal{E}(\vr, \vt | \tilde \vr, \tilde
\vt )}.
\]
}

Summing up the results of the section, we may choose $\delta > 0$ so small that relation
(\ref{w18}) takes the form
\bFormula{w22}
\intO{ \left( \frac{1}{2} \vr |\vu - \tilde \vu |^2 +
\mathcal{E}(\vr, \vt| \tilde \vr, \tilde \vt ) \right) (\tau, \cdot) } + c_1 \int_0^\tau \intO{ | \Grad \vu - \Grad \tilde \vu |^2 } \ \dt
\eF
\[
+ \int_0^\tau \intO{ \left( \frac{\vc{q} (\vt, \Grad \vt) \cdot \Grad \tilde \vt }{\vt} -
\frac{\tilde \vt}{\vt}  \frac{\vc{q} (\vt, \Grad \vt) \cdot \Grad  \vt }{\vt}
+ (\tilde \vt - \vt) \frac{ \vc{q}(\tilde \vt, \Grad \tilde \vt) \cdot \Grad \tilde \vt }{
{\tilde \vt}^2 } + \frac{\vc{q}(\tilde \vt, \Grad \tilde \vt)}{\tilde \vt} \cdot
\Grad (\vt - \tilde \vt)
\right) }
\ \dt
\]
{
\[
\leq c_2 \int_0^\tau  \intO{ \left( \frac{1}{2} \vr | \vu - \tilde
\vu |^2 + \mathcal{E}(\vr, \vt | \tilde \vr, \tilde \vt ) \right)
}  \ \dt.
\]
}

\subsubsection{Heat conductivity}

In accordance with hypothesis (\ref{m9}), we write
\[
\vc{q} (\vt, \Grad \vt) = - \kappa_0 \Grad \vt - \kappa_2 \vt^2 \Grad \vt - \kappa_3
\vt^3 \Grad \vt.
\]

We compute
\bFormula{w23}
\frac{\tilde \vt}{\vt} \frac{\kappa_0}{\vt} |\Grad \vt |^2 - \frac{\kappa_0}{\vt}
\Grad \vt \cdot \Grad \tilde \vt + \frac{\vt - \tilde \vt}{\tilde \vt} \frac{\kappa_0}{\tilde \vt}
|\Grad \tilde \vt |^2 + \frac{\kappa_0}{\tilde \vt} \Grad \tilde \vt \cdot \Grad (\tilde \vt - \vt)
\eF
\[
= \kappa_0 \Big[ \tilde \vt |\Grad \log(\vt) |^2 - \tilde \vt \Grad \log (\vt) \cdot \Grad \log(\tilde \vt) + (\vt - \tilde \vt) |\Grad \log (\tilde \vt) |^2 + \Grad \log (\tilde \vt) \cdot \Grad (\tilde \vt - \vt) \Big]
\]
\[
= \kappa_0 \Big[ \tilde \vt | \Grad \log(\vt) - \Grad \log (\tilde \vt) |^2 + (\vt - \tilde \vt) |\Grad \log (\tilde \vt) |^2 + \Grad \log (\tilde \vt) \cdot \Grad (\tilde \vt - \vt)
\]
\[
+ \tilde \vt \Grad \log (\tilde \vt) \cdot \Big( \Grad \log(\vt) - \Grad \log(\tilde \vt) \Big) \Big] =  \kappa_0 \Big[ \tilde \vt | \Grad \log(\vt) - \Grad \log (\tilde \vt) |^2
\]
\[
+ (\vt - \tilde \vt) |\Grad \log (\tilde \vt) |^2 + (\tilde \vt - \vt) \Grad \log(\tilde \vt) \cdot \Grad \log(\vt) \Big]
\]
\[
= \kappa_0 \left[ \tilde \vt | \Grad \log(\vt) - \Grad \log (\tilde \vt) |^2 + { (\vt - \tilde \vt) \Grad \log(\tilde \vt) \cdot \Grad \Big( \log(\tilde \vt ) - \log(\vt) \Big)} \right],
\]
where the second term on the right-hand side can be ``absorbed'' by the remaining integrals in (\ref{w22}).

Similarly, we get
\bFormula{w24}
\kappa_2 \tilde \vt  |\Grad \vt |^2 - \kappa_2 \vt
\Grad \vt \cdot \Grad \tilde \vt + \kappa_2 (\vt - \tilde \vt)
|\Grad \tilde \vt |^2 + \kappa_2 \tilde \vt \Grad \tilde \vt \cdot \Grad (\tilde \vt - \vt)
\eF
\[
= \kappa_2 \Big[ \tilde \vt | \Grad \vt - \Grad \tilde \vt |^2 + {(\vt - \tilde \vt) \Grad \tilde \vt \cdot \Grad (\vt - \tilde \vt)} \Big].
\]

Finally,
\bFormula{w25}
\kappa_3 \vt \tilde \vt | \Grad \vt |^2 - \kappa_3 \vt^2
\Grad \vt \cdot \Grad \tilde \vt + \kappa_3 (\vt - \tilde \vt) \tilde \vt
|\Grad \tilde \vt |^2 + \kappa_3 {\tilde \vt}^2 \Grad \tilde \vt \cdot \Grad (\tilde \vt - \vt)
\eF
\[
= \kappa_3 \Big[ \tilde \vt \vt \Grad \vt \cdot (\Grad \vt - \Grad \tilde \vt ) +
\tilde \vt \vt \Grad \vt \cdot \Grad \tilde \vt - \vt^2 \Grad \vt \cdot \Grad \tilde \vt +
(\vt - \tilde \vt) \tilde \vt | \Grad \tilde \vt |^2 + {\tilde \vt}^2 \Grad \tilde \vt \cdot \Grad (\tilde \vt - \vt) \Big]
\]
\[
= \kappa_3 \Big[ \tilde \vt \vt | \Grad \vt - \Grad \tilde \vt |^2 + 2 \vt \tilde \vt \Grad \vt \cdot \Grad \tilde \vt - \vt^2 \Grad \vt \cdot \Grad \tilde \vt - \tilde \vt^2 \Grad \vt
\cdot \Grad \tilde \vt \Big]
\]
\[
= \kappa_3 \tilde \vt \vt | \Grad \vt - \Grad \tilde \vt |^2 - \kappa_3 (\vt - \tilde \vt)^2 \Grad \vt \cdot \Grad \tilde \vt,
\]
where
\[
(\vt - \tilde \vt)^2 \Grad \vt \cdot \Grad \tilde \vt = |\Grad \tilde \vt |^2 (\vt - \tilde \vt)^2 + \Grad (\vt - \tilde \vt) \cdot \Grad \tilde \vt (\vt - \tilde \vt)^2.
\]
We conclude by observing that
\[
\Grad (\vt - \tilde \vt) \cdot \Grad \tilde \vt (\vt - \tilde \vt)^2 =
\Grad (\vt - \tilde \vt) \cdot \Grad \tilde \vt [ \vt - \tilde \vt]^2_{\rm ess} +
\Grad (\vt - \tilde \vt) \cdot \Grad \tilde \vt [ \vt - \tilde \vt]^2_{\rm res},
\]
where, furthermore,
\bFormula{w26} \intO{\left| \Grad (\vt - \tilde \vt) \cdot \Grad
\tilde \vt [ \vt - \tilde \vt]^2_{\rm res} \right|} \leq \intO{|
\Grad \tilde \vt | \left[ \delta |\Grad \vt - \Grad \tilde \vt |^2
+ c(\delta) \left[\vt - \tilde \vt \right]^4_{\rm res} \right]}
\eF
{
\[
\le \delta\|\Grad\vu-\Grad\tilde\vu\|^2_{L^2(\Omega;R^{3\times
3})}+ K(\delta,\cdot) \intO{{\cal
E}(\vr,\vt|\tilde\vr,\tilde\vt)},
\]
\[
\intO{\left| \Grad (\vt - \tilde \vt) \cdot \Grad \tilde \vt [ \vt
- \tilde \vt]^2_{\rm ess} \right|} \leq 2\overline\vt\intO{| \Grad
\tilde \vt | \left[ \delta |\Grad \vt - \Grad \tilde \vt | +
c(\delta) \left[ \vt - \tilde \vt \right]^2_{{\rm ess}} \right]}
\]
\[
\le \delta\|\Grad\vu-\Grad\tilde\vu\|^2_{L^2(\Omega;R^{3\times
3})}+ K(\delta,\cdot) \intO{{\cal E}(\vr,\vt|\tilde\vr,\tilde\vt)}
\]
for any $\delta > 0$.
}

Summing up (\ref{w23} - \ref{w26}) we may write (\ref{w22}) as
\bFormula{w27}
\intO{ \left( \frac{1}{2} \vr |\vu - \tilde \vu |^2 +
\mathcal{E}(\vr, \vt| \tilde \vr, \tilde \vt ) \right) (\tau, \cdot) } +
c_1 \int_0^\tau \intO{ | \Grad \vu - \Grad \tilde \vu |^2 } \ \dt
\eF
\[
+ c_2 \left[ \int_0^\tau \intO{ |\Grad \vt - \Grad \tilde \vt |^2
} \ \dt + \int_0^\tau \intO{ | \Grad (\log(\vt)) - \Grad
\log(\tilde \vt) |^2 } \ \dt \right]
\]
\[
\leq c_3 \int_0^\tau \intO{ \left( \frac{1}{2} \vr | \vu - \tilde
\vu |^2 + \mathcal{E}(\vr, \vt | \tilde \vr, \tilde \vt ) \right)
} \ \dt \ \mbox{for a.a} \ \tau \in (0,T),
\]
which yields the desired conclusion
\[
\vr \equiv \tilde \vr, \ \vt \equiv \tilde \vt, \ \vu \equiv \tilde \vu.
\]

We have proved Theorem \ref{Tm1}.

\section{Concluding remarks}
\label{c}

The structural restrictions introduced in Section \ref{m}, and, in particular, the presence of the radiation pressure proportional to $\vt^4$ were motivated by the \emph{existence theory} developed in \cite[Chapter 3]{FENO6}. More refined arguments could be used to show that many of these assumptions could be relaxed in the proof of weak-strong uniqueness.
In particular, the estimates based on the presence of the radiation components of the thermodynamic functions $p$, $e$, and $s$ could be performed by means of the ``dissipative'' terms on the left-hand side of (\ref{w27}) combined with some variant of Poincare's inequality.

Similar result can be obtained for other types of boundary conditions and even on unbounded spatial domains.

Last but not least, we remark that the smoothness assumptions imposed on the classical solution $\tilde \vr$, $\tilde \vt$, $\tilde \vu$ can be relaxed in terms of the integrability properties of the weak solution $\vr$, $\vt$, $\vu$.

Our final remark concerns the weak formulation of the Navier-Stokes-Fourier system introduced in Section \ref{ws}. Note that there are at least two alternative ways how to replace the entropy balance (\ref{i3}), namely by the
total energy balance
\bFormula{C1}
\partial_t \left( \frac{1}{2} \vr |\vu|^2 + \vr e(\vr, \vt) \right) +
\Div \left[ \left( \frac{1}{2} \vr |\vu|^2 + \vr e(\vr, \vt)  + p(\vr,\vt) \right) \vu \right] + \Div \vc{q} =
\Div (\tn{S} \vu ),
\eF
or by the internal energy balance
\bFormula{C2}
\partial_t (\vr e(\vr, \vt)) + \Div (\vr e(\vr, \vt) \vu ) + \Div \vc{q} = \tn{S} : \Grad \vu - p \Div \vu.
\eF
Although (\ref{C1}), (\ref{C2}) are equivalent to (\ref{i3}) for classical solutions, this is, in general, not the case
in the framework of weak solutions. As we have seen, it is precisely the entropy balance (\ref{i3}) that gives rise, in combination with (\ref{i8}), the relative entropy inequality (\ref{d8}) yielding the weak-strong uniqueness property.
This fact may be seen as another argument in favor of the weak formulation of the Navier-Stokes-Fourier system based
on (\ref{i3}), (\ref{i8}).

\def\ocirc#1{\ifmmode\setbox0=\hbox{$#1$}\dimen0=\ht0 \advance\dimen0
  by1pt\rlap{\hbox to\wd0{\hss\raise\dimen0
  \hbox{\hskip.2em$\scriptscriptstyle\circ$}\hss}}#1\else {\accent"17 #1}\fi}

\end{document}